\magnification=1200
\font\fp=cmr9
\parskip=10pt plus 1 pt

\def\L{{\cal L}}

\def\qed{\ifhmode\unskip\nobreak\fi\ifmmode\ifinner\else\hskip5 pt \fi\fi
\hbox{\hskip5 pt \vrule width4 pt  height6 pt  depth1.5 pt \hskip 1 pt }}
\overfullrule=0pt
\input amssym



\font\msbmnormal=msbm10
\font\msbmpeq=msbm7
\font\msbmmuypeq=msbm5
\newfam\numeros
\textfont\numeros=\msbmnormal
\scriptfont\numeros=\msbmpeq
\scriptscriptfont\numeros=\msbmmuypeq
\def\num{\fam\numeros\msbmnormal}
\def\N{{\num N}}
\def\R{{\num R}}
\def\C{{\num C}}
\def\Z{{\num Z}}

\def\T{{\num T}}

\def\P{{\num P}}

\def\boxit#1#2{\setbox1=\hbox{\kern#1{#2}\kern#1}%
\dimen1=\ht1 \advance\dimen1 by #1 \dimen2=\dp1 \advance\dimen2 by #1
\setbox1=\hbox{\vrule height\dimen1 depth\dimen2\box1\vrule}%
\setbox1=\vbox{\hrule\box1\hrule}%
\advance\dimen1 by .4pt \ht1=\dimen1
\advance\dimen2 by .4pt \dp1=\dimen2 \box1\relax}
\def\orlar#1{\ifmmode\boxit{3pt}{$#1$}
\else
\boxit{3pt}{#1}\fi}

\hsize 15true cm
\centerline{ \bf Measures and semi-orthogonal functions on the unit circle }

\medskip

$$\vbox{\halign
{&\quad \hfil # \hfil\quad \cr
{\it M.J.Cantero},$\;\;\;${\it M.P.Ferrer},$\;\;\;${\it L.Moral},$\;\;\;${\it L.Vel\'azquez}
\cr
\cr
\cr
{ \fp Departamento  de Matem\'atica Aplicada. Universidad de Zaragoza. Spain.}
\cr}}
$$

\medskip

{\bf Abstract}.

{\fp  The zeros of semi-orthogonal functions with respect to a probability measure $\mu$ supported
on the unit circle can be applied to obtain Szeg\H o quadrature formulas. The discrete measures
generated by these formulas weakly converge to the orthogonality measure $\mu$. In this paper we
construct families of semi-orthogonal functions with interlacing zeros, and give a representation
of the support of $\mu$ in terms of the asymptotic distribution of such zeros.}

\medskip

{\noindent \it Keywords and phrases:} { \fp Quadrature formulas, Semi-orthogonal functions,
Asymptotic distribution of zeros, Support of a measure. }

{\noindent\it (1991) AMS Mathematics Subject Classification  }: 42C05

\bigskip

{\bf \S $\;\;\;$ 1 -  Introduction. }

\medskip

Given a probability measure $\nu$ on the real line, it is well known (see for example [5] or [4])
that a sequence of orthogonal polynomials (SOP), $\bigl(p_n\bigr)_{n\in\N}$, with respect to $\nu$,
satisfies the following properties:

\item{a)} The zeros of $p_n$ are real and simple for $n\geq1$.

\item{b)} If $I$ is an interval such that ${\rm supp} \, \nu \subseteq I$, the zeros of $p_n$ lie on
the interior of $I$ for $n\geq1$.

\item{c)} The polynomials $p_n$ and $p_{n+1}$ have interlacing zeros for $n\geq1$.

\item{d)} The interpolatory quadrature formulas
$$
\int_\R F(x) \, d\nu(x) \approx \sum_{k=1}^{n} F(x_k^{(n)}) H_k^{(n)}, \qquad n\geq1,
$$
are exact for every polynomial $F$ with degree not greater than
$2n-1$ if the nodes $\bigl(x_k^{(n)}\bigr)_{k=1}^n$ are the zeros
of $p_n$ and $H_k^{(n)}=1/K_{n-1}(x_k^{(n)},x_k^{(n)})$, where
$K_n(x,y)$ is the $n$-kernel associated to
$\bigl(p_n\bigr)_{n\in\N}$, which is
$$
K_n(x,y) = {k_n \over k_{n+1}}{p_{n+1}(x) p_n(y) - p_n(x) p_{n+1}(y) \over (x-y) \int_\R p_n^2(t) \, d\nu(t)},
$$
being $k_n$ the leading coefficient of $p_n$. These quadrature
formulas imply the existence of a sequence
$\bigl(\nu_n\bigr)_{n\geq1}$ of discrete measures on the real line
given by
$$
d\nu_n(x) = \sum_{k=1}^{n} {\delta(x-x_k^{(n)}) \over K_{n-1}(x_k^{(n)},x_k^{(n)})} \, dx
$$
($\delta$ is the Dirac distribution), such that $\nu_n \buildrel
* \over \longrightarrow \nu$ when the support of $\nu$ is bounded.

Therefore, there is a close relationship between the measure $\nu$
and the location of the zeros of a SOP. In fact, properties a),
b), c), d) can be used to study the support of $\nu$ in terms of
the asymptotics of such zeros (see [5], [4]).

As for measures on the unit circle $\T$, the situation is rather
different. The associated orthogonal polynomials have their zeros
inside the open unit disk, and not on $\T$, where the measure is
supported. Moreover, these zeros are not necessarily simple. So,
it has no meaning to look for results completely analogous to the
ones above.

In the case of the real line, we can find other functions that
satisfy property a), and a weak version of b), c), d): polynomials
$q_n(x), n\geq1$, with $\deg(q_n) = n$, and orthogonal to $x^k$
for all $k<n-1$, but not to $x^{n-1}$. These are called
quasi-orthogonal polynomials (see [5], [4]), and the modified
versions of properties b), c), d) that they satisfy are:

\item{b')} If $I$ is an interval such that ${\rm supp} \, \nu \subseteq I$, it is not possible to
ensure that all the zeros of $q_n$ lie on the interior of $I$, but at most one of them lies outside.

\item{c')} It is not always true that $q_n$ and $q_{n+1}$ have interlacing zeros, but, if the convex
hull of ${\rm supp} \, \nu$ is not $\R$, the sequence
$\bigl(q_n\bigr)_{n\geq1}$ can be chosen so that all the
polynomials $q_n$ have a common zero outside this convex hull, and
the remaining zeros of $q_n$ and $q_{n+1}$ (that are inside this
convex hull) interlace.

\item{d')} The interpolatory quadrature formulas constructed by using as nodes the zeros of $q_n$
are exact for every polynomial with degree not greater than $2n-2$, instead of $2n-1$. That is,
they are exact on a $(2n-1)$-dimensional vector subspace of the space of real polynomials. Likewise,
these quadrature formulas generate a sequence of discrete measures on the real line weakly converging
to $\nu$.

As in the case of orthogonal polynomials, properties a), b'), c'),
d') allow to study the support of $\nu$ through the analysis of
the asymptotic distribution of zeros of quasi-orthogonal
polynomials.

The aim of this paper is to get on $\T$ a situation as similar as
possible to the one we have in the real line with quasi-orthogonal
polynomials: given a probability measure $\mu$ on $\T$, we look
for sequences $\bigl(F_n\bigr)_{n\geq1}$ of analytic functions on
a open subset $\Omega \supset \T$ of $\C$, such that, when
restricted to $\T$ can be viewed as real functions, and that
satisfy analogous properties to a), b'), c'), d'). That is,

\item{1.} $F_n(z)$ is real for $z\in\T$.

\item{2.} The zeros of $F_n$ are simple and lie on $\T$.

\item{3.} If $I$ is a connected subset of $\T$ such that ${\rm supp} \, \mu \subseteq I$,
at most one zero of $F_n$ lies on the closure of $\T \backslash I$.

\item{4.} The zeros of the functions $F_n$ and $F_{n+1}$ satisfy an \lq\lq interlacing property" (that
we will define exactly later on).

\item{5.} There exist interpolatory quadrature formulas
$$
\int_\T F(z) \, d\mu(z) \approx \sum_{k=1}^{n} F(z_k^{(n)}) H_k^{(n)},
\eqno(1)
$$
where $\bigl(z_k^{(n)}\bigr)_{k=1}^n$ are the zeros of $F_n$ and $H_k^{(n)}>0$, that are exact on certain
$(2n-1)$-dimensional vector subspace of the space of Laurent polynomials, and that generate a sequence
$\bigl(\mu_n\bigr)_{n\geq1}$ of discrete measures supported on $\T$ and weakly converging to $\mu$.

Notice that a SOP with respect to $\mu$ does not satisfy any of above properties. The search of these
analog of quasi-orthogonal polynomials will lead to the so-called para-orthogonal polynomials
([9], [6], [3]). Finally, with a simple modification of these polynomials we will arrive at a sequence of
Laurent polynomials called semi-orthogonal functions ([1], [2], [3]), that will provide the answer for the
requirements 1, 2, 3, 4, 5.

Once we have found these sequences of functions, the final aim of this work will be to characterize the
support of $\mu$ in terms of the asymptotic distribution of their zeros.

\bigskip

{\bf \S$\;\;\;$ 2 - Basic tools. }

\medskip

We shall denote by $\Z$ the set of all integer numbers, and
$\N=\{n\in\Z: n\geq 0\}$. Let $\P$ be the vector space of
polynomials with complex coefficients, and $\P_n$ the vector
subspace of polynomials whose degree is less than or equal to $n$.

Let $\mu$ be a probability measure supported on
$\T:=\{z\in\C:|z|=1\}$. Let us define on $\P$ the inner product
$$
\big< f, g \big> := \int_{\T} f(z) \overline{g(z)} d\mu(z),\qquad f,g\in\P.
$$
The Gram-Schmidt procedure determines the sequence of monic orthogonal polynomials (SMOP),
$\bigl(\Phi_n\bigr)_{n\in\N}$, satisfying the well known recurrence relation
$$
\eqalign{&\Phi_0(z) = 1, \cr
         &\Phi_{n+1}(z) = z\Phi_n(z) + \Phi_{n+1}(0) \Phi_n^*(z),  \qquad n\in\N,\cr}
$$
where $\Phi_n^*(z) = z^n \overline{\Phi}_n(z^{-1})$ is the reversed polynomial of $\Phi_n$.
Furthermore, for $n\geq 1$, the zeros of the polynomial $\Phi_n$ lie on $|z|<1$ and, thus,
$|\Phi_n(0)|<1$.

The sequence $\bigl(\Phi_n(0)\bigr)_{n\in\N}$ is the so called sequence of Schur parameters, and
the condition $|\Phi_n(0)|<1$, $ n\geq 1$ is in fact equivalent to the existence of an unique
probability measure for which $\bigl(\Phi_n\bigr)_{n\in\N}$ is a SMOP (see [11], [7], [5]).

If we denote $e_n:=\bigl<\Phi_n,\Phi_n\bigr> >0, \, n\geq 0$, then,
$$
\eqalign{&e_0=1, \cr
         &e_n=\prod_{k=1}^n\bigl(1-|\Phi_k(0)|^2\bigr), \qquad n\geq 1, \cr}
$$
and $\bigl(e_n\bigr)_{n\geq 0}$ constitutes a strictly decreasing
sequence. The orthonormal polynomial sequence related to above
inner product,  $\bigl(\varphi_n\bigr)_{n\geq 0}$, is
$\displaystyle \varphi_n(z)={1\over \sqrt{e_n}} \Phi_n(z)$.

The Kernel polynomials $K_n(z,y)$ are defined by
$$
K_n(z,y) := \sum_{k=0}^n\varphi_k(z)\overline{\varphi_k(y)} = \sum_{k=0}^n
{\Phi_k(x)\overline{\Phi_k(y)} \over e_k}, \qquad   n\in\N.
$$
They satisfy the reproducing property
$$
\bigl< K_n(z,y),  f(z)\bigr> = \overline{f(y)}, \qquad   f\in\P_n,
$$
as well as the Christoffel-Darboux formula
$$
e_{n+1} (1-\overline{y}z) K_n (z,y) = \Phi_{n+1}^*(z) \overline{\Phi_{n+1}^*(y)} -
\Phi_{n+1}(z)\overline{\Phi_{n+1}(y)},
\qquad   n\in\N.
$$

Now, let us denote by $\Lambda={\rm span}\{z^n: n\in\Z\}$ the vector space of Laurent polynomials
with complex coefficients, and, for $p,q\in\Z$ with $p\leq q$, let $\Lambda_{p,q}$ be the following
vector subspace of $\Lambda$
$$
\Lambda_{p,q} :=  \biggl\{\sum_{k=p}^q \alpha_k z^k: \alpha_k\in\C \biggr\}.
$$

The measure $\mu$ induces a linear functional on $\Lambda$, $\L:\Lambda\to \C$, given by
$$
\L\bigl[f\bigr] :=\int_{\T} f(z) d\mu(z),  \qquad  f\in\Lambda .
$$
An orthogonalization procedure on $\bigl(\Lambda_{-n,n},\L\bigr), \, n\in\N$, leads to
$$
{\cal B}=\{1\}\cup\biggl(\bigcup_{n\geq 1}\{f_n^{(1)}, f_n^{(2)}\}\biggr),
$$
with
$$
\left.\eqalign{
f_n^{(1)}(z)&={\overline{\alpha}_n \Phi_{2n}(z) + \alpha_n \Phi_{2n}^*(z) \over z^n}=
              {\overline{A}_n z\Phi_{2n-1}(z) + A_n \Phi_{2n-1}^*(z) \over z^n}\cr
f_n^{(2)}(z)&={\overline{\beta}_n \Phi_{2n}(z) - \beta_n \Phi_{2n}^*(z) \over i z^n}=
              {\overline{B}_n z\Phi_{2n-1}(z) - B_n \Phi_{2n-1}^*(z) \over i z^n} \cr}
\right\}
\eqno(2)
$$
where $\alpha_n,\beta_n\in\C$ are such that Re$(\alpha_n\overline{\beta}_n)\not=0$. The set
${\cal B}$ of Laurent polynomials constitutes a basis for $\Lambda$ such that
$$
\L\bigl[f_n^{(j)}(z) \, z^k\bigr] = 0,  \quad -n+1\leq k \leq n-1, \quad  j=1,2,
$$
and the matrix
$$
\bigl(\L\bigl[f_n^{(j)}(z) f_n^{(l)}(z) \bigr]\bigr)_{j,l=1,2}
$$
is positive definite for each $n\geq 1$. These functions $f_n^{(j)}, \, n\geq 1, \, j=1,2$, are called
semi-orthogonal functions (SOF) with respect to $\L$ or $\mu$ (see [1], [2]).

Notice that, when $z\in\T$, then $f_n^{(j)}(z)\in\R$. So, SOF satisfy the first property we are
looking for. To develop under which conditions they verify the rest of the requirements listed in
the introduction, we start by studying the numerator polynomials of SOF.

\medskip

{\bf Invariant para-orthogonal polynomials and quadrature formulas. }

\smallskip

{\bf Definition}.  We say that $P\in\P_n$ is an $n$-para-orthogonal polynomial ($n$-POP) with
respect to the measure $\mu$ if the following statements hold:

\item{i)} $\bigl< P, z^k\bigr> = 0, \qquad  1\leq k \leq n-1$,

\item{ii)} $\bigl< P, 1\bigr> \not=0,\qquad  \bigl< P, z^n\bigr> \not=0$.

{\bf Definition}. For $\kappa\in\C$, $\kappa\not=0$, a polynomial $P$ is called
$(\kappa)$-invariant if
$$
P^*(z) = \kappa P(z),  \qquad \forall z\in\C,
$$
where  $P^*(z) = z^{{\rm deg}(P)} \overline{P}(z^{-1})$.

{\bf Remark 1}. Notice that $|P^*(z)| = |P(z)|$ for $z\in\T$ implies that necessarily $|\kappa|=1$.

{\bf Proposition 1.}  { \sl A polynomial $P\in\P_n$ is an   invariant  $n$-POP if and only
if there exist $\alpha,\beta \in\C\setminus\{0\}$ such that $|\alpha| =|\beta| $ and
$$
P(z) = \alphaÊ\Phi_n(z) + \betaÊ\Phi_n^*(z).
$$
}
{ \bf{\sl  Proof.} } It follows immediately from the orthogonal decomposition
$$
\P_n = z \P_{n-2} \oplus  {\rm span}\{\Phi_n, \Phi_n^*\}.
$$
$\diamond\diamond$

{\bf Remark 2}. Notice that the numerator of a SOF is a 1-invariant 2$n$-POP.

POP are introduced and studied for the first time in [9]. More
recent results about them can be found in [6], [3]. Here, we will
make only a summary of the mean properties that we need for our
purposes (that is, properties relative to the zeros of POP),
referring the reader to these works for the proofs.

{\bf Proposition 2.} {\sl Invariant POP satisfy the following properties:

\item{\bf 1.}  The zeros of an invariant $n$-POP are simple and lie on $\T$ (see [9], [3]).

\item{\bf 2.}  Let $P$ be a polynomial such that $P(z)=C_1 \Phi_n(z) + C_2 \Phi_n^*(z)$
and with zeros on $\T$. Then, $P$ is an invariant  $n$-POP (see [3]).

\item{\bf 3.}  Let $P$ be an invariant $n$-POP and let $z_0\in \T$ be a zero of $P$. Then, we can
write $P(z)= c (z-z_0) K_{n-1}(z, z_0)$, with $c\not=0$.
Conversely, $c (z - z_0 ) K_{n-1}(z, z_0)$ with $c\not=0$ and $z_0\in\T$ is an  invariant $n$-POP
(see [3]).

\item{\bf 4.}  Let $P_n^{(1)}$ and $P_n^{(2)}$ be two invariant $n$-POP with respect to the same measure.
If they have a common zero, then there exists $\lambda\in\C\backslash\{0\}$ such that
$P_n^{(1)} = \lambda P_n^{(2)}$ (see [3]).

\item{\bf 5.}  The vector subspace ${\rm span}\{\Phi_n, \Phi_n^*\bigr\}$ is generated by any two
invariant $n$-POP having not a common zero (see [3]).

}

Above properties are the starting point to construct interpolatory quadrature formulas like (1) that use
as nodes the zeros of invariant POP. The results we summarize here can be found in [9]. Some extensions
of them are discussed in [6], [3].

Given a measure $\mu$ supported on $\T$, we shall call an $n$-point Szeg\"o quadrature problem, the following
problem:

{\sl Find $n$ points $\bigl(z_k\bigr)_{k=1}^n \subset \T$, $z_k\not=z_j$ if $k\not=j$ and $n$ positive
numbers $\bigl(H_k\bigr)_{k=1}^n$ such that the quadrature formula
$$
\int_{\T} F(z) \, d\mu(z) \approx \sum_{k=1}^n F(z_k) H_k,\qquad F\in C(\T),
\eqno(3)
$$

\item{(i)}  is exact for every $F\in \Lambda_{r,s}$, for some integers $r,s$,

\item{(ii)} is not exact for every $F\in \Lambda_{r',s'}$ if $\Lambda_{r',s'} \supsetneqq \Lambda_{r,s}$.

}

Requirement (ii) is just a maximality condition, and is the responsible for the impossibility of
ensuring the existence of a solution to above problem in general (see [3]). However, when $-r=s=n-1$
such a maximal solution does always exists as the following theorem asserts.

{\bf Theorem 1.}  { \sl  The $n$-point Szeg\"o quadrature problem has a solution on $\Lambda_{-n+1,n-1}$
if and only if the nodes $\bigl(z_k\bigr)_{k=1}^n$ are the zeros of an invariant $n$-POP (see [9], [6], [3]).
In such a case, the weights $H_k$ in the quadrature formula (3) are given by (see [6], [3])
$$
H_k = {1\over K_{n-1}(z_k, z_k)}, \qquad k=1, \dots, n.
$$

}

It is possible to find for the weights $H_k$ another expression in terms of POP. Let $P(z)$ be an invariant
$n$-POP whose zeros are $\bigl(z_k\bigr)_{k=1}^n$, and consider the rational function
$\displaystyle f(z) = {P(z) \over z^p}$ where $p\in\Z$ is fixed. Notice that $f \in \Lambda_{-p,n-p}$, and
has the same zeros as $P$ because $z=0$ can not be a zero of $P$. Moreover, $f$ provides a set of Laurent
polynomials $\bigl(l_k\bigr)_{k=1}^n \subset \Lambda_{-p,n-p-1}$ through the following expressions
$$
l_k(z) := {f(z) \over (z-z_k) f'(z_k)}, \qquad  k=1,\dots,n,
$$
that is a basis for $\Lambda_{-p,n-p-1}$ because $l_k(z_j) = \delta_{kj}$.
If $0 \leq p \leq n-1$, then $\Lambda_{-p,n-p-1} \subset \Lambda_{-n+1,n-1}$, and the quadrature
formula (3) that use as nodes the zeros of $P$ is exact for $l_k, \; k=1,\dots,n$.
Therefore,
$$
\int_{\T} l_k(z) \, d\mu(z) = \sum_{j=1}^n l_k(z_j) H_j = H_k, \qquad  k=1,\dots,n.
$$

Summarizing, if $P(z)$ is an invariant $n$-POP and $p$ is an integer such that $0 \leq p \leq n-1$,
the rational function $\displaystyle f(z) = {P(z) \over z^p}$ generate a solution to the
$n$-point Szeg\"o quadrature problem on $\Lambda_{-n+1,n-1}$: the nodes are the zeros of $f$ and the
weights are given in terms of $f$ by
$$
 H_k = \int_{\T} {f(z) \over (z-z_k) f'(z_k)} \, d\mu(z), \qquad  k=1,\dots,n.
\eqno(4)
$$
Of course, in spite of the $p$-dependence of $f$, the expression in the right hand side of (4) must be
independent of $p$ for $0 \leq p \leq n-1$ because the weight $H_k$ so is.

If deg$(P)$ is even, say 2$n$, we can take $p=n$, and, thus
$$
f(z) = {\overline{\alpha} \Phi_{2n}(z) + \alpha \Phi_{2n}^*(z) \over z^n}.
$$
In other words, $f(z)$ is a SOF with respect to $\mu$ and therefore real for $z\in\T$.

However, if deg$(P)=2n+1$, the choice for getting $f(z)$ real on $\T$ is $p=n+{1\over 2}$, and in this case
$$
f(z) = {\overline{\alpha} \Phi_{2n+1}(z) + \alpha \Phi_{2n+1}^*(z) \over z^n z^{1\over 2}}
$$
is not a Laurent polynomial, and, even more, is not analytic in
any open subset $\Omega\supset\T$ of $\C$. Moreover, the weights
$H_k$ can not be obtained like in (4) throughout such an $f$.

That is, if $n$ is even, the quadrature formula can be generated by a SOF, but it is not such a case,
in principle, when $n$ is odd.

\medskip

{\bf Study of the case odd ${\bf n}$}

\smallskip

Let
$$
f(z) = {\overline{\alpha} \Phi_{2n+1}(z) + \alpha \Phi_{2n+1}^*(z) \over z^n z^{1\over 2}},
$$
where the determination for $z^{1\over 2}$ is arbitrarily chosen.
We wish to find a SOF with the same zeros and that, therefore, can
be used to generate quadrature formulas.

Obviously, the numerator of this $f$ is a 1-invariant
$(2n+1)$-POP, and  $f$ vanishes at $2n+1$ different points on
$\T$. Let $w$ be one of them. Then, up to a  non zero real factor,
$$
f(z) = {z-w \over i z^{1\over 2} w^{1\over 2}} \left({w\over z}\right)^n K_{2n}(z,w),
\eqno(5)
$$
where again the determination of $w^{1\over 2}$ is arbitrarily
chosen.

Let us consider the measure $\widetilde{\mu}$, given by $d\widetilde{\mu} = |z-w|^2 d\mu$
and denote $\bigl(\psi_n\bigr)_{n\in\N}$ the associated SMOP. Then, the sequences
$\bigl(\Phi_n\bigr)_{n\in\N}$ and $\bigl(\psi_n\bigr)_{n\in\N}$ are related by means
of (see [8])
$$
(z-w) \psi_{n-1}(z) = \Phi_n(z) - {\Phi_n(w) \over K_{n-1}(w,w)} K_{n-1}(z,w).
$$
Taking reversed polynomials, we get
$$
-(z-w) \psi_{n-1}^*(z) = w\Phi_n^*(z) - {\Phi_n^*(w) \over K_{n-1}(w,w)} z K_{n-1}(z,w).
$$
Then
$$
(z-w) \bigl[\Phi_n^*(w) z  \psi_{n-1}(z) + \Phi_n(w)  \psi_{n-1}^*(z) \bigr]=
$$
$$
=\Phi_n^*(w) z \Phi_n(z) - w \Phi_n(w) \Phi_n^*(z) =
$$
$$
={e_n\over e_{n-1}} \bigl[ \Phi_{n+1}^*(w) \Phi_{n+1}(z) - \Phi_{n+1}(w) \Phi_{n+1}^*(z)\bigr]
$$
follows from the elimination of $K_{n-1}(z,w)$ on above
relationships and the recurrence formula. Thus, from (5),
$$
f(z) = {e_{2n+1}\over e_{2n}}{z-w \over i z^{1\over 2} w^{1\over 2}}
{\Phi_{2n}^*(w) z \psi_{2n-1}(z) + \Phi_{2n}(w) \psi_{2n-1}^*(z) \over (zw)^n}.
$$
Taking $z=e^{i\theta}$, $w=e^{i\beta}$, we have for $z\not=w$
$$
{z-w \over i z^{1\over 2} w^{1\over 2}} = \pm 2\sin \bigl({\theta - \beta \over 2}\bigr) \not= 0,
$$
where the sign depends on the determination of $z^{1\over 2}$ and $w^{1\over 2}$.
That is,
$$
f(z) = \pm 2 {e_{2n+1}\over e_{2n}} \sin \bigl({\theta - \beta \over 2}\bigr) {\Phi_{2n}^*(w) z
\psi_{2n-1}(z) + \Phi_{2n}(w)
\psi_{2n-1}^*(z) \over (zw)^n}.
\eqno(6)
$$
So, if we define the Laurent polynomial
$$
g(z) = \bigl({w\over z}\bigr)^n K_{2n}(z, w) \in \Lambda_{-n,n},
$$
we can write, from (5) and (6)
$$
g(z) = {e_{2n+1}\over e_{2n}} {\Phi_{2n}^*(w) z \psi_{2n-1}(z) + \Phi_{2n}(w) \psi_{2n-1}^*(z) \over (zw)^n}.
$$

Notice that $g$ has the same zeros as $f$, up to $w$. Moreover, as we can see from (2), $g$ is a SOF with
respect to $\widetilde{\mu}$. Therefore, $g(z) \in\R$ for $z\in\T$, and its numerator is a 1-invariant
$2n$-POP with respect to $\widetilde{\mu}$. Besides, $(z-w)g(z)$ is a Laurent polynomial which has the same
zeros as $f$, and its numerator is a 1-invariant $(2n+1)$-POP with respect to $\mu$.

\bigskip

{\bf \S $\;\;\;$ 3 - Semi-orthogonal functions and discrete measures. }

\medskip

For each $n\in\N$, we choose $w_n\in\T$. Then, we define a sequence $\bigl(F_n(z;w_n)\bigr)_{n\geq 1}$
in $\Lambda$ by means of
$$
\left.\eqalign{
F_{2n}(z;w_{2n}) & :=
{\Phi_{2n}^*(w_{2n})\Phi_{2n}(z) - \Phi_{2n}(w_{2n}) \Phi_{2n}^*(z)
\over i (z w_{2n})^n}, \; n\geq1 \cr
F_{2n+1}(z;w_{2n+1}) & :=
{\Phi_{2n}^*(w_{2n+1})z\psi_{2n-1}(z) + \Phi_{2n}(w_{2n+1}) \psi_{2n-1}^*(z)
\over  (z w_{2n+1})^n}, \; n\geq0 \cr}
\right\}
\eqno(7)
$$
so that $F_{2n},F_{2n+1}\in\Lambda_{-n,n}$.
This sequence is constituted by SOF: $\bigl(F_{2n}\bigr)_{n\geq 1}$ is a sequence of SOF with respect
to $\mu$, and $\bigl(F_{2n+1}\bigr)_{n\geq 0}$ is another one with  respect to the family of varying
measures $\mu^{(2n+1)}$, where
$$
d\mu^{(2n+1)}(z) = |z - w_{2n+1}|^2 d\mu(z), \qquad   z\in\T.
$$
Let us denote $\bigl(z_k^n\bigr)_{k=1}^n := \{z\in\C: F_n(z;w_n)=0\} \cup \{w_n\}$.
Now, the quadrature formulas generated by $F_{2n}(z;w_{2n})$ and $(z-w_{2n+1})F_{2n+1}(z;w_{2n+1})$
can be written together as
$$
\int_{\T} F(z) \, d\mu(z) \approx \sum_{k=1}^n F(z_k^{(n)}) H_k^{(n)},
\eqno(8)
$$
where
$$
H_k^{(n)} = {1\over K_{n-1}(z_k^{(n)}, z_k^{(n)})}.
$$
Moreover, formula (8) is exact for every $F\in\Lambda_{-n+1,n-1}$.

Above results can be read as follows. Let us consider the discrete measures on $\T$
$$
d\mu_n(e^{i\theta}) = \sum_{k=1}^n {1\over K_{n-1}(z_k^{(n)}, z_k^{(n)})}
\delta\bigl(\theta - \theta_k^{(n)}\bigr) \, d\theta, \quad  n\geq1,
$$
where $\theta_k^{(n)} = {\rm Arg}(z_k^{(n)})$, and $\bigl(z_k^{(n)}\bigr)_{k=1}^n$ are the zeros of
$$
f_n(z;w_n) = {\Phi_n^*(w_n)\Phi_n(z) - \Phi_n(w_n)\Phi_n^*(z) \over i (z w_n)^{n\over 2}},
\qquad  n\geq 1,
$$
with $\bigl(w_n\bigr)_{n\geq 1} \subset \T$.
Then, the exactness for $F\in\Lambda_{-n+1,n-1}$ means that
$$
\int_{\T} F(z) \, d\mu_n(z) = \int_{\T} F(z) \, d\mu(z)  \qquad \forall F\in\Lambda_{-n+1,n-1}.
$$
When $F=1$ we get $\displaystyle \int_{\T} d\mu_n(z) = \int_{\T} d\mu(z) =c_0 = 1$, that
is, $\bigl(\mu_n\bigr)_{n\geq 1}$ is a sequence of uniformly bounded measures.
Therefore, for every $F\in C(\T), G\in\Lambda_{-m,m}$, if $n>m$ then
$$
\eqalign{
\bigg| \int_{\T} F(z) \, d\mu(z) - \int_{\T} F(z) \, d\mu_n(z) \bigg| & \leq
\bigg| \int_{\T} (F(z)-G(z)) \, d\mu(z) \bigg| + \cr
& + \bigg| \int_{\T} (F(z)-G(z)) \, d\mu_n(z) \bigg|
\leq 2||F-G||_\infty.
}
$$
Thus, since every continuous function on $\T$ can be uniformly approximated by Laurent polynomials,
$$
\lim_n \int_{\T} F(z) \, d\mu_n(z) = \int_{\T} F(z) \, d\mu(z), \qquad \forall F\in C(\T),
$$
holds. That is, $\mu_n\to \mu$ in the $\ast$-weak topology.

Summarizing, giving an arbitrary sequence $\bigl(w_n\bigr)_{n\geq 1}\subset \T$, we can define a
sequence $\bigl(\mu_n\bigr)_{n\geq 1}$ of uniformly bounded discrete measures on $\T$ such that,
$\mu_n$ is supported on the zeros of $f_n(z;w_n)$, and $\mu_n \buildrel * \over\longrightarrow \mu$.

\bigskip

{\bf \S $\;\;\;$ 4 - Sequences with interlacing zeros. }

\medskip

Let $\mu$ be a probability measure on $\T$ and let $\bigl(\Phi_n\bigr)_{n\in\N}$ be the related
SMOP. For every sequence of complex numbers $\bigl(\alpha_n\bigr)_{n\geq 1}$, with $\alpha_n\not=0$,
we can define the functions
$$
f_n(z) = {\overline{\alpha}_n \Phi_n(z) - \alpha_n \Phi_n^*(z) \over i z^{n\over 2}},
\qquad  n\geq 1,
\eqno(9)
$$
where the determination on $z^{n\over 2}$  can be arbitrarily
chosen, because we are in\-te\-res\-ted only in the location of
their zeros. We know that $f_n$ has their zeros on $\T$,
$f_n(e^{i\theta})\in\R$, and, moreover, a discrete measure $\mu_n$
is generated such that $\mu_n \buildrel * \over\longrightarrow
\mu$.

The study of the properties of the zeros of $\bigl(f_n\bigr)_{n\geq 1}$ is equivalent to the study of
the same properties for the SOF $\bigl(F_n(z;w_n)\bigr)_{n\geq 1}$ given in (7).

By using the recurrence relation for $\bigl(\Phi_n\bigr)_{n\in\N}$ in (9), we have
$$
f_n(z) = {e_n\over e_{n-1}}
{\overline{\beta}_n z \Phi_{n-1}(z) - \beta_n \Phi_{n-1}^*(z) \over i z^{n\over 2}},
\qquad  n\geq 1,
\eqno(10)
$$
where $\displaystyle \beta_n={e_{n-1}\over e_n} (\alpha_n - \Phi_n(0)\overline{\alpha}_n)$,
$\,\displaystyle \overline{\beta}_n={e_{n-1}\over e_n} (\overline{\alpha}_n - \overline{\Phi_n(0)} \alpha_n)$,
and thus
$$
\alpha_n = \beta_n + \Phi_n(0) \overline{\beta}_n, \qquad  n\geq 1.
\eqno(11)
$$
Without loss of generality, we can assume that $\bigl(\alpha_n\bigr)_{n\geq 1}$ is given by means
of a polynomial sequence $\bigl(p_n\bigr)_{n\geq 1}$ in the following way
$$
\alpha_n:= w^{-{m\over 2}} p_n(w),
$$
where $w\in\T$ is arbitrary, and $m = m(n) =$ deg$(p_n)$, which is
chosen so that $m(n+1)=m(n)+1$. Under these assumptions, we can
write
$$
\beta_n = w^{-m\over 2}{e_{n-1}\over e_n} \bigl[p_n(w) - \Phi_n(0) p_n^*(w)\bigr].
$$

Let $\omega_0\in\R$ be fixed. We can write $\displaystyle
\T=\{e^{i\theta}: \theta\in[\omega_0,\omega_0+2\pi)\}$, and, thus,
$f_n(e^{i\theta})$ is a $C^{\infty}$ real function defined on
$[\omega_0, \omega_0+2\pi)$, where the non integer powers of $z$
and $w$ are taken so that ${\rm arg}(z^{1\over 2}), {\rm
arg}(w^{1\over 2}) \in \bigl[{\omega_0\over 2}, {\omega_0\over
2}+\pi\bigr)$. For two numbers $\zeta_1$, $\zeta_2\in\T$, let us
denote $\theta_j = {\rm arg}(\zeta_j) \in  [\omega_0,
\omega_0+2\pi)$. Then, we can establish an order relation by
$$
\zeta_1 < \zeta_2 \iff \theta_1 < \theta_2.
$$
Let $\bigl(\zeta_j^{(n)}\bigr)_{j=1}^n$ and $\bigl(\zeta_j^{(n+1)}\bigr)_{j=1}^{n+1}$ be the zeros of
$f_n$ and $f_{n+1}$ respectively, which are ordered as above. Then, we wish to determine a suitable sequence
$\bigl(\alpha_n\bigr)_{n\geq 1}$ such that, for each $n\geq 1$,
$$
\zeta_j^{(n+1)} < \zeta_j^{(n)} < \zeta_{j+1}^{(n+1)}, \qquad \;\;\; j=1,\dots, n,
$$
holds when $\zeta_j^{(k)} \not= e^{i\omega_0}$, $k=n,n+1$. In short, in such a case we say that $f_n$
and $f_{n+1}$ have interlacing zeros. It is a well known result, derived from Sturm's theorem, that
the interlacing condition is equivalent to
$$
{\rm sgn}\bigl[ \bigl({{\rm d} \over {\rm d}\theta} f_{n+1}(\zeta) \bigl) f_n(\zeta)\bigr] =
{\rm sgn}\bigl[ i \zeta f_{n+1}^{'}(\zeta) f_n(\zeta)\bigr] =
{\rm constant}
$$
for all $\zeta\in \bigl(\zeta_j^{(n+1)}\bigr)_{j=1}^{n+1}$, $\zeta\not=e^{i\omega_0}$.

Let us consider for $f_n$ and $f_{n+1}$ the expressions (9) and (10), respectively,
$$
\eqalign{
&f_n(z) = {\overline{\alpha}_n \Phi_n(z) - \alpha_n \Phi_n^*(z) \over i z^{n\over 2}},\cr
&f_{n+1}(z) = {e_{n+1}\over e_n}  {\overline{\beta}_{n+1} z \Phi_n(z) - \beta_{n+1} \Phi_n^*(z)
\over i z^{n+1\over 2}}.   \cr}
$$
Since $f_{n+1}(\zeta)=0$, we have
$\displaystyle {\beta_{n+1}\over \overline{\beta}_{n+1}} = {\zeta\Phi_n(\zeta)\over \Phi_n^*(z)}$.
Thus
$$
\beta_{n+1} = k_{n+1} \zeta^{-n+1\over 2}Ê\Phi_n(\zeta), \qquad
\overline{\beta}_{n+1} = k_{n+1} \zeta^{-n-1\over 2}Ê\Phi_n^*(\zeta),
\eqno(12)
$$
with $k_{n+1}=k_{n+1}(\zeta)\in\R\setminus \{0\}$. For $f_{n+1}$, and using Christoffel-Darboux formula,
we get
$$
f_{n+1}(z) = k_{n+1} e_{n+1} {z-\zeta\over i\zeta} \left({\zeta\over z}\right)^{n\over 2} K_n(z,\zeta).
$$
Hence,
$$
i\zeta f_{n+1}^{'}(\zeta) = k_{n+1} e_{n+1} K_n(\zeta,\zeta).
$$
On the other hand, from (12),
$$
f_n(\zeta) = {1\over k_{n+1}} {\overline{\alpha}_n \beta_{n+1} \zeta^{-{1\over 2}} -
\alpha_n \overline{\beta}_{n+1} \zeta^{1\over 2} \over i}
$$
follows. That is,
$$
i\zeta f_{n+1}^{'}(\zeta) f_n(\zeta) = 2 e_{n+1} K_n(\zeta,\zeta) \,
{\overline{\alpha}_n \beta_{n+1} \zeta^{-{1\over 2}} -
\alpha_n \overline{\beta}_{n+1} \zeta^{1\over 2} \over 2i}.
$$
The last factor in the right hand side suggests the choice
$$
w=e^{i\omega_0}, \qquad \beta_{n+1} = w^{1\over 2}\alpha_n,
\eqno(13)
$$
that leads to an expression with constant sign
$$
{\overline{\alpha}_n \beta_{n+1} \zeta^{-{1\over 2}} -
\alpha_n \overline{\beta}_{n+1} \zeta^{1\over 2} \over 2i} =
|\alpha_n|^2 \sin \bigl({\omega_0-\theta \over 2 }\bigr) < 0,
$$
where $\theta={\rm arg} (\zeta)\in (\omega_0, \omega_0+2\pi)$.

Now, from (13) and (11) it follows that
$$
p_{n+1}(w) = w p_n(w) + \Phi_{n+1}(0) p_n^*(w), \qquad  n\geq 1.
\eqno(14)
$$
This equation has two independent solutions (see [7], [5]), that
can be chosen as the sequences $\bigl(\Phi_n(w)\bigr)_{n\geq 1}$
and $\bigl(i\Omega_n(w)\bigr)_{n\geq 1}$, where $\Omega_n$ is the
$n$-th second kind polynomial. Thus, there exist two fixed
polynomials $A$, $B$, and an integer $k\geq 0$, such that
$$
\left\{ \eqalign{
          &p_n(w) = A(w) \Phi_n(w) + B(w) \Omega_n(w),  \cr
          &A^{*_k}(w) = A(w), \quad  B^{*_k}(w) = - B(w),  \cr}\right.
\eqno (15)
$$
where $P^{*{_k}}(z) = z^k \overline{P}(z^{-1})$, deg$(P)\leq k$. The conditions (14) and (15) are
equivalents (see [10]).

Summarizing, we can deduce the following result.

{\bf Theorem 2.}  { \sl  Let $\bigl(p_n\bigr)_{n\in\N}$ be a polynomial sequence satisfying (14), or
equivalently (15). Let $w=e^{i\omega_0}$ be fixed. Then, for the sequence of functions
$\bigl(f_n\bigr)_{n\geq 1}$ defined by
$$
f_n(z) = {\overline{p_n(w)} \Phi_n(z) - p_n(w) \Phi_n^*(z) \over i z^{n\over 2}}, \qquad  n\geq1,
\eqno(16)
$$
$f_n$ and $f_{n+1}$ have interlacing zeros in $(\omega_0,\omega_0+2\pi)$.}

As two particular cases of (16), we have
$$
\eqalign{
f_n^{(1)} (z;w) &:={\Phi_n^*(w) \Phi_n(z) - \Phi_n(w) \Phi_n^*(z) \over i (zw)^{n\over 2}},
\qquad  n\geq 1,  \cr
f_n^{(2)} (z;w) &:={\Omega_n^*(w) \Phi_n(z) + \Omega_n(w) \Phi_n^*(z) \over  (zw)^{n\over 2}},
\qquad  n\geq 1.  \cr}
\eqno (17)
$$
From (15) and (16),
$$
\eqalign{
f_n(z) = &{A(w) \over w^{k\over 2}} f_n^{(1)} (z;w) - {B(w) \over w^{k\over 2}} f_n^{(2)} (z;w) = \cr
       = & A_1 f_n^{(1)} (z;w) + A_2 f_n^{(2)} (z;w), \cr}
$$
(where $\displaystyle A_1={A(w) \over w^{k\over 2}}, A_2={-{ B(w) \over w^{k\over2}}}$) follows.
Moreover, $A_1$ and $A_2$ are real numbers, independent from $n$, with $(A_1,A_2)\not=(0,0)$.

{\bf Remark 3}. Conversely, given $ \displaystyle f_n(z) = {\overline{\alpha}_n \Phi_n(z) -
\alpha_n \Phi_n^*(z) \over i z^{n\over 2}}, \, \alpha_n\not=0$, there exist two real numbers
$A_1(n)$, $A_2(n)$ with $(A_1(n), A_2(n))\not=(0,0)$ such that
$$
 f_n(z)=A_1(n) f_n^{(1)} (z;w) + A_2(n) f_n^{(2)} (z;w).
$$
In fact, let us denote
$$
\eqalign{
&P_n(z) = \overline{\alpha}_n \Phi_n(z) - \alpha_n \Phi_n^*(z),  \cr
&P_n^{(1)}(z)  = w^{-{n\over 2}}\bigl[\Phi_n^*(w) \Phi_n(z) - \Phi_n(w) \Phi_n^*(z) \bigr],\cr
&P_n^{(2)}(z)  = w^{-{n\over 2}}\bigl[\Omega_n^*(w) \Phi_n(z) + \Omega_n(w) \Phi_n^*(z) \bigr].\cr}
$$
These polynomials are 1-invariant $n$-POP. Moreover,
$$
P_n \in {\rm span} \{\Phi_n,\Phi_n^*\} = {\rm span} \{P_n^{(1)}, P_n^{(2)}\}.
$$
Thus, there exists $(A_1(n), A_2(n))\not=(0,0)$ such that
$$
P_n(z) = A_1(n) P_n^{(1)}(z) + {A_2(n)\over i} P_n^{(2)}(z)
$$
holds. Thus,
$$
f_n(z) = A_1(n) f_n^{(1)}(z;w) + A_2(n) f_n^{(2)}(z;w).
$$
By taking $z=w$ and $z\in\T$, it follows that $A_1(n),A_2(n)\in\R$.

The choice (13), which guarantees interlacing property, corresponds to the case when $A_1(n)$, $A_2(n)$
are independent of $n$. Of course, (13) is not the unique option for interlacing zeros.

{\bf Remark 4}. Notice that the functions $f_n^{(1)}(z;w)$ have a common zero at $\zeta_1^{(n)}=w$.
In this case, the interlacing property is strictly  verified in the interval $(\omega_0,\omega_0+2\pi)$.
By using the well-known identity
$$
\Omega_n^*(z) \Phi_n(z) + \Omega_n(z) \Phi_n^*(z) = 2 e_n z^n
\eqno(18)
$$
(see [7]), we obtain that $f_n^{(2)}(w;w) = 2 e_n$. This implies that, for the remaining cases, i.e.,
when $B(w)\not=0$ in (15), $w$ is not a zero of $f_n$. Thus, the interlacing property holds in
$[\omega_0, \omega_0+2\pi)$.

Finally, we are able to give other kind of interlacing properties.

{\bf Proposition 3.}  { \sl For each $n\geq 1$, the functions $f_n^{(1)}(z;w)$ and $f_n^{(2)}(z;w)$
have interlacing zeros in $[\omega_0, \omega_0+2\pi)$.}

{\bf{\sl Proof.}} If is enough to prove that the expression
$$
i \zeta f_n^{{(1)^{'}}}(\zeta;w) f_n^{(2)}(\zeta;w)
$$
has the same sign for every zero $\zeta$ of $f_n^{(1)}(z;w)$.
If $f_n^{(1)}(\zeta;w)=0$, then
$$
{\Phi_n(\zeta) \over \Phi_n(w)} = {\Phi_n^*(\zeta) \over \Phi_n^*(w)}.
$$
Using (18) we get
$$
\eqalign{
f_n^{(2)}(\zeta;w) &= {\Omega_n^*(w) \Phi_n(w) + \Omega_n(w) \Phi_n^*(w) \over (\zeta w)^{n\over 2}}
{\Phi_n^*(\zeta) \over \Phi_n^*(w)} = \cr
&= 2 e_n {\Phi_n^*(\zeta) \over \Phi_n^*(w)} \left({w\over \zeta}\right)^{n\over 2}. \cr}
$$
On the other hand,
$$
\eqalign{
f_n^{(1)}(z;w) &={\Phi_n^*(w) \over \Phi_n^*(\zeta)} \left({\zeta\over w}\right)^{n\over 2}
{\Phi_n^*(\zeta)\Phi_n(z) - \Phi_n(\zeta) \Phi_n^*(z) \over i (z\zeta)^{n\over 2}}=\cr
&={\Phi_n^*(w) \over \Phi_n^*(\zeta)} \left({\zeta\over w}\right)^{n\over 2}f_n^{(1)}(z;\zeta)=\cr
&={\Phi_n^*(w) \over \Phi_n^*(\zeta)}  {\zeta\over (wz)^{n\over 2}} e_n
K_{n-1}(z, \zeta){z-\zeta\over i\zeta}. \cr}
$$
Thus,
$$
i \zeta f_n^{{(1)^{'}}}(\zeta;w) = {\Phi_n^*(w) \over \Phi_n^*(\zeta)}
\left({\zeta\over w}\right)^{n\over 2} e_n K_{n-1}(\zeta,\zeta).
$$
Finally,
$$
i \zeta f_n^{{(1)^{'}}}(\zeta;w) f_n^{(2)}(\zeta;w) = 2e_n^2 K_{n-1}(\zeta,\zeta) >0
$$
for each zero $\zeta$ of $f_n^{(1)}(z;w)$. $\diamond\diamond$

{\bf Corollary 1.}  { \sl Let $A_j$, $B_j \in \R, \; j=1,2$, and
$$
\eqalign{
 f_n(z) &= A_1 f_n^{(1)}(z;w) + A_2 f_n^{(2)}(z;w),  \cr
 g_n(z) &= B_1 f_n^{(1)}(z;w) + B_2 f_n^{(2)}(z;w),  \cr}
$$
where $\left|\matrix{A_1&A_2\cr B_1&B_2\cr}\right|\not= 0$. Then, $f_n$ and $g_n$ have interlacing
zeros in $[\omega_0, \omega_0+2\pi) $.}

{\bf{\sl Proof.}} Notice that $f_n$ and $g_n$ are linearly independents. Then, from Proposition 2, they
have not a common zero.
Let $\zeta_0$ be such that  $f_n(\zeta_0)=0$. Then, up to a non-zero real factor,
$f_n(z)=f_n^{(1)}(z;\zeta_0)$ holds. Thus, again from Proposition 2, we have that there exist
$C_1$, $C_2 \in \R$, with $C_2\not=0$, such
that $g_n(z) =  C_1 f_n^{(1)}(z;\zeta_0) + C_2 f_n^{(2)}(z;\zeta_0) $ . Hence, by using
Proposition 3, we have for each zero $\zeta$ of $f_n$ that
$$
\eqalign{
i \zeta f_n^{'}(\zeta) g_n(\zeta) &= i \zeta f_n^{{(1)^{'}}}(\zeta;\zeta_0)
C_2f_n^{(2)}(\zeta;\zeta_0) =\cr
&= 2 e_n^2 K_{n-1}(\zeta,\zeta) C_2,\cr}
$$
which keeps constant sign. $\diamond\diamond$

\bigskip

{\bf \S $\;\;\;$ 6 - Semi-orthogonal functions and support of the orthogonality measure. }

\medskip

We have find that any sequence of functions given in Theorem 2 satisfies properties 1, 2, 4, 5
listed in {\bf \S 1}. It remains to reach the property 3, as well as to develop the relation
between ${\rm supp}\,\mu$ and the asymptotic distribution of zeros of such kind of functions.
The starting point for this relation are the quadrature formulas described in {\bf \S 2}, that are
constructed by using as nodes the zeros of functions with the general form
$$
f_n(z) = {\overline{\alpha}_n \Phi_n(z) - \alpha_n \Phi_n^*(z) \over i z^{n\over 2}},
\qquad \alpha_n \neq 0.
\eqno(19)
$$
As we said in {\bf \S 3}, these quadrature formulas provide discrete measures $\mu_n$ weakly converging
to the orthogonality measure $\mu$, and supported on the zeros of the functions $f_n$. As a consequence,
given an arbitrary sequence $\bigl(f_n\bigr)_{n \geq 1}$ of functions with the form (19), a point of
${\rm supp}\,\mu$ must be a limit point of the set of zeros of all the functions $f_n$, or a zero of
infinitely many functions $f_n$. Thus, denoting by $A'$ the derived set of $A$, we have the following
proposition.

{\bf Proposition 4.} { \sl Let $\bigl(f_n\bigr)_{n \geq 1}$ be any sequence of functions with
the form (19), and let ${\cal Z}$ and ${\cal X}$ be the set of zeros of all the functions $f_n$,
and the set of complex numbers that are zeros of infinitely many functions $f_n$, respectively. Then,
$$
{\rm supp}\,\mu \subseteq {\cal Z}' \cup {\cal X}.
$$
}

Since ${\cal Z}'' \subset {\cal Z}'$ and ${\cal X}' \subset {\cal
Z}'$, a corollary follows immediately.

{\bf Corollary 2.} { \sl Under the assumptions of Proposition 4,
$$
({\rm supp}\,\mu)' \subseteq {\cal Z}'.
$$
}

Corollary 2 says that a limit point of ${\rm supp}\,\mu$ must be a
limit point of ${\cal Z}$ too. From Proposition 4, the rest of the
points of ${\rm supp}\,\mu$, that is, the isolated points, could
be limit points of ${\cal Z}$, or zeros of infinitely many
functions $f_n$. Now, we are going to demonstrate that, indeed, an
isolated point $z_0$ of ${\rm supp}\,\mu$ that is not a limit
point of ${\cal Z}$, must be a zero of all the functions $f_n$
except, at most, finitely many of them. If it were not the case,
then there would be a subsequence $\bigl(f_{n_k}\bigr)_{k\in\N}$
such that $f_{n_k}(z_0)\neq0$ for $k\in\N$. Moreover, we are
supposing that $z_0$ is an isolated point of ${\rm supp}\,\mu$
that is not a limit point of ${\cal Z}$, so, there must be a
neighborhood ${\cal U}$ of $z_0$ where none of the functions
$f_{n_k}$ can vanish, and such that ${\cal U}\cap{\rm
supp}\,\mu=\{z_0\}$. Therefore, $\int_{{\cal U}\cap\T} d\mu =
\mu(\{z_0\}) \neq 0$. Now, let us consider the subsequence of
measures $\bigl(\mu_{n_k}\bigr)_{k\in\N}$ corresponding to the
subsequence of functions $\bigl(f_{n_k}\bigr)_{k\in\N}$. Since
$\bigl(\mu_n\bigr)_{n\in\N}$ is uniformly bounded,
$\bigl(\mu_{n_k}\bigr)_{k\in\N}$ so is, and Helly's theorem
implies the existence of a subsubsequence
$\bigl(\mu_{n_{k_j}}\bigr)_{j\in\N}$ which is convergent to $\mu$
a.e.. For this subsubsequence it must be
$$
\lim_j \int_{{\cal U}\cap\T} d\mu_{n_{k_j}} = \int_{{\cal U}\cap\T} d\mu.
$$
But, since the functions $f_{n_k}$ have no zeros in ${\cal U}$, it is
${\cal U} \cap {\rm supp}\,\mu_{n_{k_j}} = \emptyset$, and, hence,
$\int_{{\cal U}\cap\T} d\mu_{n_{k_j}} = 0$, which is in contradiction with
$\int_{{\cal U}\cap\T} d\mu \neq 0$.

Then, we can state the following result.

{\bf Proposition 5.} { \sl Under the assumptions of Proposition 4,
$$
{\rm supp}\,\mu \subseteq {\cal Z}' \cup {\cal \widetilde X},
$$
where ${\cal \widetilde X}$ is the set of complex numbers that are zeros of all the functions $f_n$
except, at most, finitely many of them.
}

To go beyond these results we need a property that is quasi-reciprocal of previous one, and,
therefore, leads to an inclusion in opposite sense. Proposition 5 says that if a connected
subset of $\T$ is \lq\lq eventually free of zeros", then its interior is
\lq\lq out of the support of $\mu$''. On the contrary, the following proposition implies that
if a connected subset of $\T$ is \lq\lq out of the support of $\mu$'', then its closure contains
at most one zero.

{\bf Proposition 6.}  { \sl Let $\bigl(f_n\bigr)_{n \geq 1}$ be any sequence of functions with the
form (19), and let $z_1, z_2$ be two zeros of $f_n$. Then, in any of the two connected components of
$\T \backslash \{ z_1, z_2 \}$ there exits at least one point of ${\rm supp}\,\mu$, as well as one
zero of $f_m$ for $m>n$.}

{ \bf{\sl  Proof.} } Since we are just interested in the location of zeros we will only work with
the polynomials $P_n(z) = \overline{\alpha}_n \Phi_n(z) - \alpha_n \Phi_n^*(z)$ that have the same
zeros as $f_n$. Notice first that
$\displaystyle{zP_n(z) \over (z-z_1)(z-z_2)} \in {\rm span} \{ 1,z,z^2, \dots , z^{n-1} \}$, and,
therefore,
$$
\int_{\T} {zP_n(z) \over (z-z_1)(z-z_2)} \overline{P}_n(z^{-1}) \, d\mu(z) = 0,
$$
due to the orthogonality properties of $\Phi_n$ and $\Phi_n^*$.
If $\theta_j={\rm Arg}(z_j)$, above result reads
$$
\int_0^{2\pi} {|P_n(e^{i\theta})|^2 \over
\sin ({\theta-\theta_1 \over 2}) \sin ({\theta-\theta_2 \over 2})}
\, d\mu(e^{i\theta}) = 0.
$$
But the expression inside the integral has now constant and opposite sign in each connected
component of $\T \backslash \{ z_1, z_2 \}$, except for the zeros of $f_n$, for which is null
(including $z_1$ and $z_2$). Thus, taking into account that ${\rm supp}\,\mu$ is an infinite set,
if the integral has to vanish, then there must be points of ${\rm supp}\,\mu$ in both components.

On the other hand, we have that
$\displaystyle {zP_n(z) \over (z-z_1)(z-z_2)} \overline{P}_n(z^{-1}) \in \Lambda_{-n+1,n-1}$,
so, the quadrature formulas described in {\bf \S 3} are exact over this Laurent polynomial if we
choose as nodes the zeros $z_k^{(m)}$ of $f_m$ for $m>n$. Hence, if
$\theta_k^{(m)} = {\rm Arg}(z_k^{(m)})$, it must be
$$
\sum_{k=1}^{m} {|P_n(e^{i\theta_k^{(m)}})|^2 \over
\sin ({\theta_k^{(m)}-\theta1 \over 2}) \sin ({\theta_k^{(m)}-\theta_2 \over 2})} \, H_k^{(m)} = 0,
\qquad m > n.
$$
If all the zeros $z_k^{(m)}$ were located in only one of the connected components of
$\T \backslash \{ z_1, z_2 \}$, all the non-null terms in this sum would have the same sign.
Since there must be non-null terms because the number of zeros of $P_n$ is $n<m$, this is
in contradiction with the fact that above sum must vanish. $\diamond\diamond$

Notice that this proposition has as a direct consequence the property 3 listed in {\bf \S 1} that we
were looking for. That is,

{\bf Corollary 3.} { \sl If $I$ is a connected subset of $\T$ such that ${\rm supp}\,\mu \subseteq I$,
at most one zero of $f_n$ lies on the closure of $\T \backslash I$.}

To get more accurate results about ${\rm supp}\,\mu$ we need to restrict our attention to those
sequences of functions with interlacing zeros. Then, we will suppose that $\bigl(f_n(z)\bigr)_{n\geq1}$
is one of the sequences described in Theorem 2, that is,
$$
f_n(z) = A_1 f_n^{(1)} (z;w) + A_2 f_n^{(2)} (z;w), \qquad  (A_1,A_2)\neq(0,0),
\eqno(20)
$$
where $f_n^{(i)}, \, i=1,2$, are given in (17). In that case, $f_n$ and $f_{n+1}$ can have only one
common zero at the point $w$, due to the interlacing of zeros in $\T\backslash\{w\}$. Thus,
${\cal \widetilde X}\subseteq\{w\}$. In fact, from Remark 4, we can see that
${\cal \widetilde X}=\{w\}$ only when $A_2=0$. Thus, we have reached the following result.

{\bf Proposition 7.} { \sl Let $\bigl(f_n\bigr)_{n \geq 1}$ be any sequence of functions with the form
(20). Then, if $A_2=0$,
$$
{\rm supp}\,\mu \subseteq {\cal Z}' \cup \{w\},
$$
and, if $A_2\neq0$,
$$
{\rm supp}\,\mu \subseteq {\cal Z}'.
$$
}

Now, we will work only with the functions $f_n^{(1)}(z;w)$, that is, those functions
$f_n(z)$ with a common zero at $z=w$. We will suppose for ${\rm supp}\,\mu$ a non-trivial situation, what
means that ${\rm supp}\,\mu \neq \T$. The strategy will be to choose a point $w$ in each connected component
of $\T \backslash {\rm supp}\,\mu$, and consider the corresponding functions $f_n^{(1)}(z;w)$. So that,
the representation of ${\rm supp}\,\mu$ will be given throughout the zeros of several sequences
$\bigr(f_n^{(1)}(z;w)\bigr)_{n\geq1}$, one for each connected component of $\T \backslash {\rm supp}\,\mu$.

{\bf Theorem 3.} { \sl Suppose that ${\rm supp}\,\mu \neq \T$ and let
$$
\T \backslash {\rm supp}\,\mu = \bigcup_{k \in {\cal A}} C_k
$$
be the decomposition of $\T \backslash {\rm supp}\,\mu$ into connected components. If $w_k \in C_k$
for every $k \in {\cal A}$, then,
$$
{\rm supp}\,\mu = \bigcap_{k \in {\cal A}} {\cal Z}(w_k)',
$$
where ${\cal Z}(w)$ is the set of zeros of all the functions $f_n^{(1)}(z;w)$.
}

{ \bf{\sl  Proof.} } From Proposition 7 it is clear that, since $w_k \notin {\rm supp}\,\mu$, it has to be
${\rm supp}\,\mu \subseteq \bigcap_{k \in {\cal A}} {\cal Z}(w_k)'$. On the other hand, Proposition 6
implies that $w_k$ must be the unique zero of $f_n^{(1)}(z;w_k)$ in $C_k$. Therefore,
${\cal Z}(w_k)' \cap C_k = \emptyset$, and, thus,
$\bigcap_{k \in {\cal A}} {\cal Z}(w_k)' \subseteq {\rm supp}\,\mu$. $\diamond\diamond$

\bigskip

{\bf Acknowledgements}.- This research was partially supported by Universidad de Zaragoza PAI UZ-97-CIE-10
and Direcci\'on General de Ense\~nanza superior of Spain (DGES) Project PB98-1615.

\eject

{\bf References}

\medskip

\item{[1]} ALFARO, M.; CANTERO, M.J.; MORAL, L. \lq\lq
Semi-orthogonal functions and orthogonal polynomials on the unit
circle", J. Comput. Appl. Math. {\bf 99} (1998) 3-14.

\item{[2]} CANTERO, M.J. \lq\lq Polinomios ortogonales sobre la
circunferencia unidad. Modificaciones de los par\'ametros de
Schur", Doctoral Dissertation, Universidad de Zaragoza, 1997.

\item{[3]}  CANTERO, M.J.; FERRER, M.P.; MORAL, L. ; VELAZQUEZ, L.
\lq\lq Funciones semiortogonales y formulas de cuadratura", Actes
des $VI^{\grave emes}$ Journ\'ees de Math\'ematiques Apliqu\'ees
et de Statistiques (Jaca, 1999) pp. 141--150, Publications de
l'Universit\'e de Pau, 2001.

\item{[4]}  CHIHARA, T.S. \lq\lq An introduction to orthogonal
polynomials", Gordon and Breach, New York, 1978.

\item{[5]}  FREUD, G. \lq\lq  Orthogonal polynomials", Pergamon
Press, Oxford, 1971.

\item{[6]} BULTHEEL, A; GONZ\'ALEZ-VERA, P; HENDRIKSEN, E;
NJ\"ASTAD, O. \lq\lq  Orthogonal rational functions", University
Press, Cambridge, 1999.

\item{[7]}  GERONIMUS, Ya.L. \lq\lq  Orthogonal polynomials",
Consultants Bureau, New York, 1961.

\item{[8]} GODOY, E; MARCELLAN, F. \lq\lq Orthogonal polynomials
and rational modifications of measures", Canad. J. Math. {\bf 45}
(1993) 930-943.

\item{[9]} JONES,  W.B.; NJ\"ASTAD, O; THRON, W.I. \lq\lq Moment
theory, orthogonal polynomials, quadrature, and continued
fractions associated with the unit circle", Bull. London Math.
Soc. {\bf 21} (1989) 113-152.

\item{[10]} PEHERSTORFER, F; STEINBAUER,  R. \lq\lq
Characterization of general orthogonal polynomials with respect to
a functional", J. Comp. Appl. Math. {\bf 65 } (1995) 339-355.

\item{[11]} SZEG\H O, G. \lq\lq Orthogonal Polynomials", Amer.
Math. Soc. Colloq. Publ. 23, Providence, Rhode Island, 1975
(Fourth Edition).

\bye
\end